\documentclass[11pt,a4paper]{article}
\usepackage[utf8]{inputenc}
\usepackage[T1]{fontenc}
\usepackage{lmodern}
\usepackage[english]{babel}
\usepackage[margin=2.5cm]{geometry} 
\usepackage{graphicx} 
\usepackage{booktabs} 
\usepackage{natbib} 
\usepackage{amsmath} 
\usepackage{amssymb} 
\usepackage{setspace} 
\usepackage{caption} 
\usepackage[colorlinks=true, citecolor=blue, linkcolor=blue, urlcolor=blue]{hyperref} 
\usepackage{amsthm}

\newtheorem{definition}{Definition}[section]

\newtheorem{conjecture}{Conjecture}[section]
\newcommand{\bed}{\begin{definition}}
\newcommand{\eed}{\end{definition}}

\newtheorem{theorem}{Theorem}

\newcommand{\eps}{\epsilon}

\newcommand{\bitem}{\begin{itemize}}
\newcommand{\eitem}{\end{itemize}}

\newcommand{\goto}{\rightarrow}

\newcommand{\beqn}{\begin{equation}}
\newcommand{\eeqn}{\end{equation}}
\newcommand{\balign}{\begin{align}}
\newcommand{\ealign}{\end{align}}

\usepackage{amssymb}
\newcommand{\beq}{\begin{equation}}
\newcommand{\eeq}{\end{equation}}

\usepackage{hyperref}
\hypersetup{
    colorlinks=true,
    linkcolor=blue,
    filecolor=magenta,      
    urlcolor=cyan,
}

\setcitestyle{round}

\let\oldfootnote\footnote
\renewcommand{\footnote}{\fontsize{9}{11}\selectfont\oldfootnote}

\title{Applications of Royen's Gaussian Correlation Inequality in high-dimensional testing}
\author{
    Jiashun Jin \\
   School of Statistics and Data Science\\
     Southeast University \\
    \and
    Zheng Tracy Ke \\
    Department of Statistics\\
     Harvard University \\
    \and
    Jingming Wang\thanks{Corresponding author: pdw9qv@virginia.edu} \\
    Department of Statistics\\
    University of Virginia
}
\date{} 

\begin{document}
\maketitle

\begin{abstract}
Royen's Gaussian Correlation Inequality (GCI) is a fundamental result in statistics and probability, but unfortunately, at least by now, the result is not very well known in the statistical community. In this paper, we demonstrate the potential of GCI in high-dimensional testing problems that involve multivariate Gaussian  variables.  Leveraging the GCI, we establish sharp Type I error bounds under arbitrary Gaussian covariance structures.  We also formulate a conjecture   and relate it to the well-known Spherical Cap Covering Conjecture in discrete geometry and Weak Simplex Conjecture in information theory.  
\end{abstract}

\section{Introduction} 
The Gaussian Correlation Inequality (GCI) is one of the most celebrated results in multivariate analysis.  Consider an $n$-dimensional zero-mean Gaussian vector $Z \sim N(0, \Sigma)$, where 
$\Sigma$ is an $n \times n$ correlation matrix (so all entries on the main diagonal of $\Sigma$ are $1$).  The GCI  claims that for any (centrally) symmetric and convex 
sets $G$ and $F$ in $\mathbb{R}^n$  
\[
\mathbb{P}(Z \in  G \cap F) \geq \mathbb{P}(Z \in G) \mathbb{P}(Z \in F). 
\] 
The GCI (in various forms) was conjectured to be true for more than 50  years,  but only by 2014, a definite and general proof was delivered by Thomas Royen  
\citep{royen2014proof}.   
Royen's proof was later popularized and verified globally by \citep{latala2017royen}, 
cementing the GCI as a fundamental pillar for non-asymptotic inference in high-dimensional statistics.  
Unfortunately, though  $10$ years have passed,  the powerful result of 
GCI still  remains largely unknown in our community. In this paper, we consider several high-dimensional test settings and show that the GCI can be used to obtain tight Type I error bounds that are valid for 
arbitrary correlation structures. We hope that our paper can help popularize the GCI within our community.

Suppose we observe an $n$-dimensional Gaussian vector 
\begin{equation} \label{model1a} 
X = \mu + Z, \qquad \mbox{where $\mu = \mathbb{E}[X]$ and $Z \sim N(0, \Sigma)$ for a correlation matrix $\Sigma$}.  
\end{equation} 
This naturally arises in genetics and genomics (e.g., \cite{IFPCA}),   where we have 
two groups of patients (e.g., normal and diseased) measured on a total of $n$ features (e.g., genes, SNPs), and $X_i$ 
is the two-sample $t$-test statistic associated with feature $i$, $1 \leq i \leq n$.  Since each $t$-test statistic 
is approximately $N(0,1)$, we assume $Z \sim N(0, \Sigma)$ and $\Sigma$ has a 
unit diagonal (so it is a correlation matrix). For theoretical purposes, such an idealization is 
conventional. 

We are interested in the problem of global testing where the goal is to test 
\begin{equation} \label{model1b} 
H_0^{(n)}: \mu = 0 \qquad vs. \qquad H_1^{(n)}: \mu \neq 0. 
\end{equation} 
When $\Sigma = I_n$ (the $n \times n$ identity matrix), the problem was well-studied (see \cite{DJ15} for a review), focusing on the so-called rare/weak signal setting. In such a setting, for a fixed $0 < a < 1$ and $\eps_n, \tau_n > 0$, 
we assume the mean vector $\mu = (\mu_1, \mu_2, \ldots, \mu_n)'$ is generated by 
\begin{equation} \label{model1c} 
\mu_i = \left\{ 
\begin{array}{ll} 
\tau_n,  &\qquad \mbox{with probability $(1 - a)   \eps_n $},  \\ 
- \tau_n, &\qquad \mbox{with probability $a  \eps_n$},  \\ 
0, &\qquad \mbox{with probability $(1 - \eps_n)$},  
\end{array} 
\right. 
\end{equation} 
Moreover, for fixed parameters $0 < \beta < 1$ and $r > 0$,  we calibrate $\eps_n$ and $\tau_n$  by 
\begin{equation} \label{model1d} 
\eps_n =  n^{-\beta}, \qquad 
\tau_n  = \left\{
\begin{array}{cc}
n^{r}, &\qquad 0 < \beta < 1/2, \\ 
\sqrt{2 r \log (n)}, &\qquad 1/2 < \beta < 1.  
\end{array} 
\right. 
\end{equation} 
Introduce the detection boundary by 
\[
\rho(\beta) = \left\{
\begin{array}{lc} 
\beta - 1/2, &\qquad 1/2 < \beta \leq 3/4, \\
(1 - \sqrt{1 - \beta})^2, &\qquad 3/4 < \beta < 1,  
\end{array} 
\right. \qquad \mbox{when $1/2 < \beta < 1$}, 
\] 
and 
\[
\rho_a^*(\beta) = \left\{
\begin{array}{cc} 
 \beta - 1/2,  &\qquad a \neq 1/2,  \\ 
(\beta - 1/2)/ 2, &\qquad a = 1/2, \\
\end{array} 
\right. 
\qquad \mbox{when $0 < \beta < 1/2$}.   
\] 
Consider the two-sided max-type test statistic 
\begin{equation} \label{max-test}  
T_n = \max\{|X_1|, |X_2|, \ldots, |X_n|\}, 
\end{equation} 
and two moment-based test statistics 
\begin{equation} \label{moment-test} 
S_n = n^{-1/2} |\sum_{i = 1}^n X_i|, \qquad  
Q_n = (2n)^{-1/2} [\sum_{i = 1}^n X_i^2 - n].  
\end{equation}
The max-type test statistic is closely related to the popular minimum p-value (MinP) test. 
To use each of these for our testing problem (\ref{model1a})-(\ref{model1d}), 
we reject the null if and only if the test statistic exceeds a threshold $t  > 0$. 
With a properly selected threshold for each test,  we have the following claims (see \cite{DJ04, DJ15} for details):  
\begin{itemize} 
\item When $1/2 < \beta < 1$ and $r > \rho(\beta)$,  the Higher Criticism 
test (e.g., \cite{DJ04}) has asymptotically full power. If we further restrict $\beta$ to 
the range of $3/4 < \beta < 1$, then the max-type test $T_n$ also has asymptotically full power. 
\item When $0 < \beta < 1/2$ and $r > \rho_a^*(\beta)$, 
the test $S_n$ has asymptotically full power if $a \neq 1/2$, but may lose power if $a = 1/2$. When $a = 1/2$,  the test $Q_n$ has asymptotically full power.  
\end{itemize} 
 

Note that the above results are for $\Sigma = I_n$. In more general cases,   the behavior of these test statistics is largely unknown. This motivates  a revisit to these test 
statistics. In this paper, we are particularly interested in the two-sided max-type statistic $T_n$ and its one-sided variant 
\[
T_n^* = \max\{X_1, X_2, \ldots, X_n\},  
\] 
which is preferred over $T_n$ in applications when most of the nonzero $\mu_i$ are positive. 
We  also consider a hybrid test $H_n$, where for thresholds $t_1, t_2, t_3 > 0$, we reject the null if and only if 
\[
\mbox{$T_n \geq t_1$ or $S_n \geq  t_2$  or $Q_n \geq t_3$}.  
\]
Note that each of these tests $T_n, S_n, Q_n$ is powerful in some regimes but may lose power in others, 
so we may want to combine them; this explains why the hybrid test $H_n$ is desirable.    
Finally, note that the Higher Criticism test is much more complicated than $T_n, S_n$ and $Q_n$, 
and is hard to analyze in the current setting where $\Sigma$ can be any correlation matrix.   
For these reasons,  we omit the Higher Criticism test for our study below.

\noindent {\bf Remark 1}. Suppose $H_0^{(n)}$ holds, i.e., $\mu = 0$. Note that  \cite{sidak1967rectangular} and \cite{khatri1967certain}  proved that  for any $a_1, \ldots, a_n\geq 0$, 
\begin{align}\label{Sidak-Khatri}
\mathbb P \big( |X_1|\leq a_1, \ldots, |X_n|\leq a_n\big) \geq \prod_{i=1}^n \mathbb P (|X_i|\leq a_i)\,. 
\end{align}
If we let $R$ be the hypercube 
\begin{equation} \label{DefineR} 
R = \cap_{i = 1}^n R_i,
\end{equation} 
where $R_i = \{x \in \mathbb{R}^n: |x_i| \leq a_i\}$ is the $i$-th coordinate slab 
 bounded by  two hyperplanes $\{x = (x_1, \ldots, x_n):  x_i =a_i\}$ and 
 $\{x = (x_1, \ldots, x_n): x_i = -a_i\}$, 
their result can be restated as that the probability that $X$ lies in the hypercube is at least the product 
of the probabilities that $X$ lies in each of the $n$ slabs. 
  In comparison, the GCI applies to much broader settings, where if the event is the intersection of $K$ ($K$ can be any integer)  
  centrally symmetric and convex sets: 
\[
G = \cap_{k = 1}^K G_k, 
\] 
then 
\[
\mathbb{P}(X \in G) \geq \Pi_{i = 1}^K \mathbb{P}(X \in G_k). 
\] 
For the three testing problems we consider above, we can use (\ref{Sidak-Khatri}) to address the first one. 
However, for more general settings such as the hybrid test, (\ref{Sidak-Khatri}) does not apply, 
and we must use Royen's GCI.

\section{Main results: tight Type I error control} \label{sec:main} 
For each of the three testing approaches, $T_n, H_n$ and $T_n^*$, 
we discuss how to find a tight control on the Type I errors. 
Consider $T_n$ first. To use the statistic for our testing problem, we pick a threshold $t > 0$ and reject $H_0^{(n)}$ if and only if 
\begin{equation} \label{MTT} 
T_n >  t. 
\end{equation} 
Recall that in our setting $X = \mu + Z$ and $Z \sim N(0, \Sigma)$. Let $F_n(x; \Sigma)$ be the cumulative distribution function of $\max_{1 \leq i \leq n} \{|Z_i|\}$.  
The Type I error for the max-type test in (\ref{MTT}) is  
\begin{equation} \label{MTT1} 
1 - \mathbb{P}(T_n \leq t | \mu = 0) = 1 -  \mathbb{P}(\max_{1 \leq i \leq n} \{|Z_i|\}  \leq t) \equiv 1 - F_n(t; \Sigma). 
\end{equation} 
Now, fix $0 < \alpha < 1$ and  we wish to control the Type I error at level $\alpha$.  
When $\Sigma = I_n$,  
\[
\mbox{$F_n(x; I_n) = 0$ if $x \leq 0$ and $F_n(x; I_n) = [\mathbb{P}(|z| \leq x)]^n$ if $x > 0$}, \qquad \mbox{where $z \sim N(0,1)$}.   
\]  
This motivates to take the threshold $t = t_n$ in a way such that  
\begin{align}\label{Definetn}
\mathbb{P}(|z| \leq t_n) = (1 - \alpha)^{1/n}. 
\end{align}

Applying the GCI, we obtain the following result, which provides tight Type I error control for $T_n$ under arbitrary $\Sigma$. 
\begin{theorem} \label{thm:two-sided} 
For any correlation matrix $\Sigma$ and any  $x$, $F_n(x; \Sigma) \geq F_n(x; I_n)$.  
Moreover, fix $0 < \alpha < 1$ and consider the max-type test (\ref{MTT}) with $t = t_n$ in (\ref{Definetn}). 
For any correlation matrix $\Sigma$, the Type I error is no more than $\alpha$. Also,  when $\Sigma = I_n$, the Type I error is exactly $\alpha$.  
\end{theorem}  

\noindent {\bf Remark 2}. Let $R$ be the hypercube defined in (\ref{DefineR})  with $a_1 = a_2= \ldots =a_n = x$.   Then we have 
\[
F(x, \Sigma) = \mathbb{P}(Z \in R).   
\] 
In this special case,  Royen's GCI reduces to (\ref{Sidak-Khatri}), so  we can also prove Theorem \ref{thm:two-sided} using the \v{S}id\'ak--Khatri inequality.  But as mentioned in Remark 1,  in more complicated  settings,  
(\ref{Sidak-Khatri}) may not apply, and we must use Royen's GCI.

Consider $H_n$ next. Recall that the test rejects $H_0^{(n)}$ if and only if 
$T_n > t_1$   or  $S_n  >  t_2$ or $Q_n >  t_3$,  
where $t_1, t_2, t_3 > 0$ are thresholds. Since $T_n, S_n, Q_n$ are dependent in a complicated way, it is hard to evaluate the Type I error precisely.  
Introduce events 
\[ 
A = \cap_{i = 1}^n \{|X_i| \leq t_1\}, \qquad   B = \{|\sum_{i = 1}^n X_i|  \leq \sqrt{n} t_2\},  \qquad  C =  \{\sum_{i =1}^n X_i^2  \leq \sqrt{2n} t_3 + n\}, 
\] 
where each event defines a (centrally) symmetric convex set in $\mathbb{R}^n$. 
Under $H_0^{(n)}$, by GCI, 
\begin{equation} \label{GIC2} 
\mathbb{P}(A \cap B \cap C) \geq \mathbb{P}(T_n \leq t_1) \mathbb{P}(S_n \leq t_2) \mathbb{P}(Q_n \leq t_3) \geq 
[\mathbb{P}(|z| \leq t_1)]^n \mathbb{P}(|S_n| \leq t_2) \mathbb{P}(Q_n \leq t_3).  
\end{equation} 
We now analyze $S_n$ and $Q_n$.  
Let $\lambda_1 \geq  \lambda_2 \geq  \ldots \geq  \lambda_n \geq 0$ be the eigenvalues of $\Sigma$ and assume 
\begin{equation} \label{Q-condition} 
\max_{1 \leq i \leq n} \lambda_i^2 / (\sum_{i = 1}^n \lambda_i^2) \goto 0.  
\end{equation} 
Introduce $a_n = {n^{-1} {\bf 1}_n' \Sigma {\bf 1}_n}$ and $b_n = n^{-1} \mathrm{trace}(\Sigma^2)$.  
Similarly as before, let $t_n$ be defined by (\ref{Definetn}). For constants $x_0, y_0 > 0$, we set  
\begin{equation} \label{hybrid-threshold} 
t_1 = t_n, \qquad t_2 = \sqrt{a_n} x_0, \qquad t_3 = \sqrt{b_n} y_0.   
\end{equation} 
\begin{theorem} \label{thm:hybrid} 
Consider the hybrid test $H_n$ where we reject the null when $T_n > t_1$ or $S_n > t_2$ or $Q_n > t_3$ with $t_1, t_2, t_3$ given in (\ref{hybrid-threshold}).   
Then the Type I error of the test is no more than
\[
1 - (1 - \alpha) \mathbb{P}(S_n \leq \sqrt{a_n} x_0) \mathbb{P}(Q_n \leq \sqrt{b_n} y_0) = 1 - (1 -\alpha) \mathbb{P}(|z| \leq x_0) \mathbb{P}(z \leq y_0) +o(1), 
\] 
which is valid for all $\Sigma$ such that ${\bf 1}_n' \Sigma {\bf 1}_n\neq 0$.
\end{theorem} 
Note that in the special case of $\Sigma = I_n$, $a_n = b_n =1$. 
In the general case, first,  Theorem \ref{thm:hybrid} provides an interesting theoretical guidance. 
Second, since in many practical situations $X$ is the vector of summary $z$-score, 
we may be able to estimate $a_n$ and $b_n$ consistently. Therefore, Theorem \ref{thm:hybrid} 
provides an interesting way to control the Type I error.

Consider the one-sided max-type test  $T_n^*$.  
Let $F_n^*(x; \Sigma)$ be the CDF of $\max\{Z_1, Z_2, \ldots, Z_n\}$, where we recall 
$Z \sim N(0, \Sigma)$.  
In view of Theorem \ref{thm:two-sided}, one may tend to believe that,  for any correlation matrix $\Sigma$, 
\begin{equation} \label{cdf1} 
{F}^*_n(x; \Sigma) \geq {F}^*_n(x; I_n). 
\end{equation} 
This is actually true if we restrict $\Sigma$ to the class of correlation matrices that are (entry-wise) non-negative. 
In fact, in this case,  (\ref{cdf1}) can be deduced from the well-known Slepian's lemma \citep{Slepian1962}. 
However, if we allow $\Sigma$ to have negative entries, then (\ref{cdf1}) may cease to 
hold.  Introduce an $n \times n$ matrix $\Sigma^*$ by 
$\Sigma^* = \frac{1}{n-1}(n I_n - {\bf 1}_n {\bf 1}_n')$  (note that $\mathrm{rank}(\Sigma^*) = n-1$).  
By basic algebra, there exists $n \times (n-1)$ matrix $ V^* = [v_1^*, v_2^*, \ldots, v_n^*]'$ such that 
\[
\Sigma^* = V^* (V^*)', \qquad \mbox{and so $\|v_i^*\| = 1$ and $(v_i^*, v_j^*) = - 1/(n-1)$ for all $i \neq j$}.   
\] 
Therefore, $v_1^*, v_2^*, \ldots, v_n^*$ are the vertices of a regular simplex inscribed in the unit sphere of  
$\mathbb{R}^{n-1}$.  Suppose $u$ is uniformly distributed over the unit sphere. For any 
$n$ points $v_1, v_2, \ldots, v_n$ on the unit sphere, introduce the random variables 
\[
U_n = \max\{u' v_1, u'v_2, \ldots, u' v_n\}, \qquad \mbox{and} \qquad U_n^* = \max\{u' v_1^*, u' v_2^*, \ldots, u' v_n^*\}. 
\] 
\begin{conjecture} \label{conj1}  
For any $x$, $\mathbb{P}(U_n \leq x) \geq \mathbb{P}(U_n^* \leq x)$. Also, the equality holds for all $x$ if and only if 
$v_1, v_2, \ldots, v_n$ are the vertices of a regular simplex inscribed in the unit sphere in $\mathbb{R}^{n-1}$.  
\end{conjecture} 
\begin{conjecture}  \label{conj2} 
For any $n \times n$ correlation matrix $\Sigma$, $F_n^*(x; \Sigma) \geq F_n^*(x; \Sigma^*)$.  
\end{conjecture} 
When $n \leq  4$,  the conjectures can be proved (see Section~\ref{sec:appen}  for details). When $n \geq 5$, both 
conjectures stand open, and it is highly non-trivial to prove or disprove them. 
In fact, the conjectures are  closely related to three famous conjectures: 
the {\it Spherical Cap Covering (SCC)} conjecture (\cite{balitskiy2017optimality}) and the 
 {\it Simplex Mean Width (SMW)} conjecture in 
convex geometry (\cite{litvak2017around}),  and the {\it Weak Simplex (WS)} conjecture (\cite{balitskiy2017optimality}) in information theory.  
With our notations above, the SMW and the WS conjectures assert that 
 \[
 \mathbb{E}[U_n] \leq \mathbb{E}[U_n^*], \qquad \mbox{and} \qquad \mathbb{E}[e^{\lambda RU_n}] \leq \mathbb{E}[e^{\lambda RU_n^*}],  
 \]  
respectively, where $R^2 \sim \chi_{n-1}^2$ is independent of $u$, and $\lambda>0$ is a fixed constant.   The SCC conjecture asserts that, for $n$ sphere caps with the same (spherical) radius on the unit sphere centered at $v_1, v_2, \ldots, v_n$, the total area 
covered by these sphere caps is maximized when $v_1, v_2, \ldots, v_n$ are the vertices of a regular simplex inscribed the unit sphere. It is seen that the $5$ conjectures are related as follows (where $\Longrightarrow$ means implies)  
\[
\mbox{SCC  $\Longrightarrow$ Conjecture \ref{conj1} $\Longrightarrow$ all other three (SMW, WS, and Conjecture \ref{conj2})}. 
\] 

Back to the one-sided max-type  statistic $T_n^*$. Fix $0 < \alpha < 1$ and suppose $t_n^* > 0$ satisfies  
\begin{equation} \label{one-sided-test} 
F_n^*(t_n^*; \Sigma^*) = (1 - \alpha).  
\end{equation} 
Since $\Sigma^*$ is completely explicit, the left hand side is known. 
Suppose we reject $H_0^{(n)}$ if and only if $T_n^* \geq t_n^*$.  
Once Conjecture 2.2 is proved, we have the following theorem, which extends Slepian's
lemma to more general cases. 
\begin{theorem} \label{thm:Slepian} 
For any correlation matrix $\Sigma$, the Type I error of the one-sided max-type test in (\ref{one-sided-test}) is no more than $\alpha$. Also, when $\Sigma = \Sigma^*$, 
the Type I error of the test is exactly $\alpha$. 
\end{theorem}

\section{Extension to multiple testing: tight FWER control} \label{sec:FWER} 
The global testing problem we considered above is closely related to the problem of  large-scale multiple testing.   
Recall that $X = \mu + Z$ with $Z \sim N(0, \Sigma)$ where $\Sigma$ is a correlation matrix. 
The goal of multiple testing is to 
simultaneously test 
\begin{equation} \label{null1} 
H_{0i}: \mu_i = 0 \qquad vs.  \qquad H_{1i}: \mu_i \neq 0,   
\end{equation} 
for each $1 \leq i \leq n$.  The Hard Thresholding (HT) test is a popular approach. For a threshold $t > 0$, the approach 
rejects $H_{0i}$ if and only if $|X_i| > t$.   In multiple testing, it is frequently desirable to control the Family-wise Error Rate (FWER) (i.e., the probability of at least one false positive; a false positive is a falsely rejected null).  The FWER associated with the HT test is 
\[
\mbox{FWER}  = 1 - \mathbb{P}(\max_{1 \leq i \leq n: \mu_i = 0} |X_i| \leq t). 
\] 
Since $\Sigma$ is unknown, how to control the FWER remains a difficult problem. 
But with GCI, we can resolve the problem as follows.  Suppose out of the $n$ null hypotheses in (\ref{null1}), 
a total of $n_0$ of them are true. Note that in many large-scale multiple testing settings, $n_0 / n \approx 1$.   It follows that in the special case of $\Sigma = I_n$, 
\[
\mbox{FWER} = 1 -  \mathbb{P}(\max_{1 \leq i \leq n: \mu_i = 0} |X_i| \leq t)=  1 -  [\mathbb{P}(|z| \leq t)]^{n_0}, \qquad \mbox{where $z \sim N(0,1)$}. 
\] 
Now, fix $0 < \alpha < 1$ and suppose we want to control the FWER at level $\alpha$. Set the threshold by $t = t_n$, where $t_n$ is determined by 
\begin{equation} \label{HT2} 
\mathbb{P}(|z| \leq t_n) =  (1 - \alpha)^{1/n}. 
\end{equation} 
Using Royen's GCI, we have the following theorem. 
\begin{theorem} \label{thm:FWER} 
Fix $0 < \alpha < 1$ and consider the HT test for (\ref{null1}) with $t_n$ in (\ref{HT2}).    
For any correlation matrix $\Sigma$, the FWER is no more than $\alpha$.  
In the special case of $\Sigma = I_n$, the FWER is $1 - (1 - \alpha)^{n_0 /n}$, which is 
nearly $\alpha$ when $n_0 / n \approx 1$.  
\end{theorem} 
As in Remark 2,  in this particular case, Royen's GCI reduces to  (\ref{Sidak-Khatri}), 
so Theorem~\ref{thm:FWER} can also be viewed as a statistical application of the \v{S}id\'ak--Khatri inequality.

\section{Numerical illustration}

\begin{figure}[th]
\centering
\begin{minipage}[c]{0.50\textwidth}
    \centering
    \includegraphics[width=0.75\linewidth]{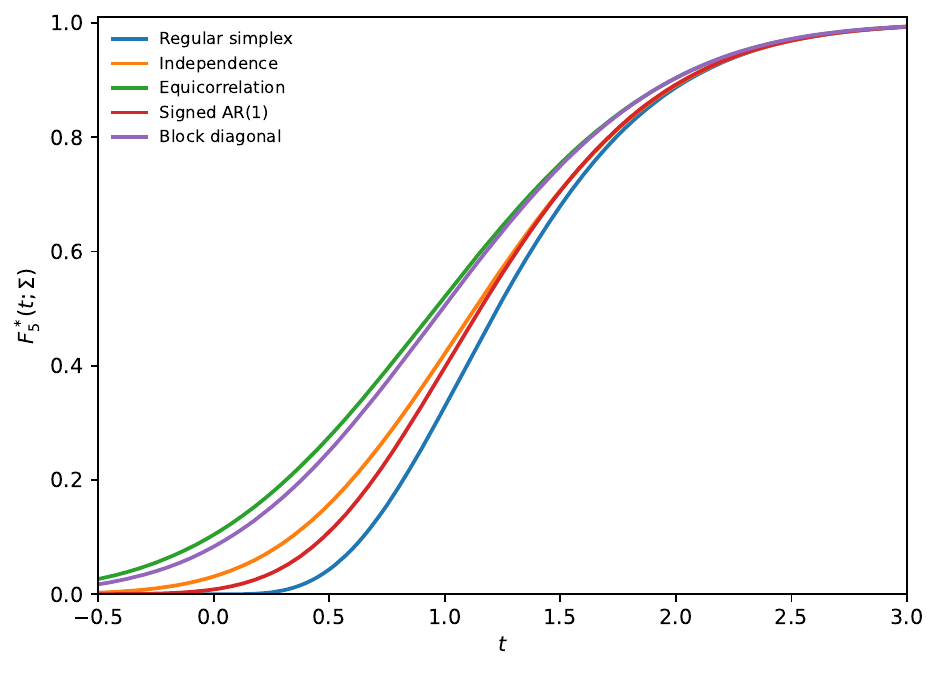}

    \small (a) CDFs of the one-sided max statistic
\end{minipage}
\hfill
\begin{minipage}[c]{0.48\textwidth}
    \centering
    \small
    \begin{tabular}{lcc}
        \toprule
        $\Sigma$ & Emp. Type I & GCI bound \\
        \midrule
        Independence    & 0.0467 & 0.0499 \\
        Equicorrelation & 0.0294 & 0.0497 \\
        Signed AR(1)    & 0.0467 & 0.0502 \\
        Block diagonal  & 0.0353 & 0.0500 \\
        \bottomrule
    \end{tabular}

    \vspace{0.6em}
    \small (b) Type I error of the hybrid test
\end{minipage}

\caption{
Numerical illustrations for the one-sided max statistic and the hybrid
test.
(a) Empirical CDFs of $\max_{1\leq i\leq5}Z_i$ under five representative
 correlation structures.
(b) Empirical Type I error of the hybrid test and the corresponding
GCI-based upper bound under four representative dependence structures.
}
\label{fig:numerical}
\end{figure}

We first provide a numerical illustration of Conjecture~\ref{conj2}.
We fix $n=5$ and consider five correlation structures:
(1) the regular-simplex correlation matrix $\Sigma^*$;
(2) independence, $\Sigma=I_5$;
(3) equicorrelation with $\rho=0.3$;
(4) signed AR(1), with $\Sigma_{ij}=(-0.5)^{|i-j|}$; and
(5) a block-diagonal structure with block sizes $2$ and $3$, within-block
correlation $0.5$, and zero between-block correlations.
For each setting, we generate $10^6$ independent Gaussian vectors and
evaluate the empirical CDF of $\max_{1\leq i\leq5}Z_i$ over
$t\in[-0.5,3]$. As shown in Figure~\ref{fig:numerical}(a), the
regular-simplex CDF lies below those under all the other structures, supporting Conjecture~\ref{conj2}.

We next illustrate the sharpness of the GCI-based Type I error bound for the hybrid
test. We take $n=50$ and consider the correlation structures (2)--(5)
above, with two equal-sized blocks for the block-diagonal case.
For each setting, the thresholds are calibrated so that $
\mathbb P(A)=\mathbb P(B)=\mathbb P(C)=0.95^{1/3} $
up to Monte Carlo error, and the GCI-based Type I error bound is
approximately $0.05$. We then estimate the empirical Type I error using an
independent  sample of size $10^6$.
As shown in Figure~\ref{fig:numerical}(b), the empirical Type I error is
below the GCI bound in all settings considered. The bound is particularly
sharp under independence and signed AR(1), while being more conservative
under positive equicorrelation and block dependence. These results illustrate that the GCI provides not only a uniform  guarantee for the hybrid test under arbitrary $\Sigma$, but can also yield a relatively sharp bound under some representative dependence structures.

\section{Discussion}

Royen's GCI is a powerful tool in probability and statistics. In this paper, we use it 
to  study the Type I error control in several high-dimensional hypothesis testing problems. While the corresponding results are well understood in the independent setting, much less is known for arbitrary correlation structures.  Leveraging Royen's GCI, we establish sharp Type I error bounds under general Gaussian correlations.

\bibliographystyle{plainnat} 
\bibliography{network} 

\appendix
\section{Appendix}\label{sec:appen} 

\noindent{\bf Further discussion of the conjectures}. We discuss how Conjecture~\ref{conj1} implies Conjecture~\ref{conj2}, SMW and WS conjectures.
Let
$
\mathcal G = \{\Sigma\in \mathbb R^{n,n}: \text{positive semidefinite,} {\rm diag}(\Sigma) = I_n \}.
$
We first claim that for any $\Sigma \in \mathcal G$ that is positive definite, we can always find a singular $\widetilde \Sigma \in \mathcal G$ such that 
$F_n^*(t; \Sigma) \geq F_n^*(t;\widetilde \Sigma)$, uniformly for $t \in\mathbb R$. Fix any   positive definite $\Sigma\in \mathcal G$. It must hold that $\max_{i\neq j}{|\Sigma_{ij}|}<1$; Otherwise, if without loss of generality $\Sigma_{12} = 1$, we have $y' \Sigma y =0$ by setting $y = (1, -1, 0, \ldots, 0)$. This implies $\Sigma$ is singular. Similarly we can claim for $\Sigma_{12} =-1$.

Given this $\Sigma$, we then introduce 
$
\widetilde \Sigma(s) = \Sigma - s({\bf 1}_n{\bf 1}_n' - I_n)\,. 
$
Notably, when $s>0$ is small such that $0<s<\lambda_{\min}(\Sigma)/(n-1)$, 
\[
\lambda_{\min}(\widetilde \Sigma(s)) \geq \lambda_{\min}(\Sigma) - s(n-1)>0\,.
\]
Thus, $\widetilde \Sigma(s)$ is positive definite and always with unit diagonal. More importantly, 
$
\widetilde \Sigma(s) \leq \Sigma
$
entrywisely. By Slepian's lemma (\cite{Slepian1962}), it yields that $F_n^*(t; \Sigma) \geq F_n^*(t; \widetilde \Sigma(s))$, uniformly for $t \in\mathbb R$.
Moreover, 
we can gradually increase $s$ until $\widetilde \Sigma(s)$ becomes singular, and define 
\[
s^* = \inf\{s>0: \lambda_{\min} (\widetilde \Sigma(s))=0\}
\]
Let $s_0 = {\bf 1}_n'\Sigma {\bf 1}_n/[n(n-1)]$. By some computations, we obtain 
$
{\bf 1}'_n\widetilde \Sigma(s_0) {\bf 1}_n = 0\,,
$
which leads to $s^*\leq s_0$. It follows from  eigenvalue continuity 
that $s^*$ exists and $\widetilde \Sigma(s^*)$ is positive semi-definite with unit diagonal. Thus,  $\widetilde \Sigma(s^*)\in \mathcal G$, and by Slepian's lemma  (\cite{Slepian1962}), 
$
F_n^*(t; \Sigma) \geq F_n^*(t; \widetilde \Sigma(s^*))$,  uniformly for $t \in\mathbb R$. 
Hence, it suffices to study 
$
 {\rm argmin}_{\Sigma\in \mathcal G^*}\mathbb{P}_{X\sim N(0, \Sigma)} (\max_{i} X_i \leq t) 
$, 
where 
$
\mathcal G^* = \{\Sigma\in \mathbb R^{n,n}: \text{ positive semidefinite, singular and} \,\, {\rm diag}(\Sigma) = I_n \}\,. 
$

Note that for any $\Sigma\in \mathcal G^*$, it holds that ${\rm rank}(\Sigma)\leq n-1$. We then represent $X_1, X_2, \ldots, X_n$ using $(n-1)$-dim Gaussian vector $Z= (Z_1, \ldots Z_{n-1})'\sim N(0, I_{n-1})$. There exists $V = (v_1, \ldots, v_n)' \in R^{n, n-1}$ such that $VV' = \Sigma$, $\|v_i \|=1$ and 
$
X_i = v_i' Z$, for $i = 1, \ldots , n$.
Using  polar decomposition of a Gaussian vector,  we have  $Z = R u$ where $R^2 \sim \chi^2_{n-1}$ independent of $u \sim Unif(\mathcal S^{n-2})$. Thus, 
$
\max_i X_i = R \cdot \max_i v_i'u
$
and
\begin{align}\label{eq:obj}
\mathbb P \big( \max_i X_i\geq t \big)  = \int \mathbb P \big(\max_i v_i'u \geq t/r \big) dF_R(r)
\end{align}
Let $\mathsf S = \{(v_1, \ldots, v_n)\in \mathcal (S^{n-2})^n: v_i'v_j =-1/(n-1), \text{ for any} \, i\neq j \}$  denote the collection of all regular simplices inscribed in a unit ball.
If we can claim that, 
\begin{align}\label{conjecture:0721}
\text{ $P \big(\max_{i\in [1,n]} v_i'u\geq t ) $  is maximized at $W = (w_1, \ldots, w_n)\in \mathsf S$ for every $t\in(0, 1]$}
\end{align}
 it follows immediately that for any $t>0$
 \begin{align}
 \mathbb P_{X\sim N(0, \Sigma)} (\max_i X_i \leq t)& = 1-  \mathbb P_{X\sim N(0, \Sigma)} (\max_i X_i \geq t)  = 1 -\int \mathbb P \big(\max_i v_i'u \geq t/r \big) dF_R(r) \notag\\
 &\geq 1 -\int \mathbb P \big(\max_i w_i'u \geq t/r \big) dF_R(r) =  \mathbb P_{X\sim N(0, \Sigma^*)} (\max_i X_i \leq t),
 \end{align}
  where $\Sigma^* = \frac{n}{n-1}I_n - \frac{1}{n-1} {\bf 1}_n{\bf 1}_n'$. For $t\leq 0$, since $\Sigma^*{\bf 1}_n= 0$, we easily get $\sum_{i=1}^n X_i = 0$ a.s. and 
  \begin{align*}
   \mathbb P_{X\sim N(0, \Sigma^*)} (\max_i X_i \leq t) 
   &=      \mathbb P_{X\sim N(0, \Sigma^*)} (X_i =0, i\in [1, n])=0\leq   \mathbb P_{X\sim N(0, \Sigma)} (\max_i X_i \leq t)
  \end{align*}
 Therefore, $F_n^*(t; \Sigma) \geq F_n^*(t;\Sigma^*)$ uniformly for $t\in \mathbb R$. As a result, it suffices to prove (\ref{conjecture:0721}). Remarkably, this can be shown rigorously for  $n= 3, 4$ by relating $P \big(\max_{i\in [1,n]} v_i'u\geq t ) $ to  the surface area of the union of spherical caps. Specifically, define the spherical cap on $\mathcal S^{n-2}$ centered at $v_i$ as
$
 C(v_i,t)=\{v\in\mathcal S^{n-2}: v_i' v\geq t\}$ for $t>0$.
 Then, we obtain 
\[
P \big(\max_{i\in [1,n]} v_i'u\geq t )  =[{\rm Vol}(\mathcal S^{n-2})]^{-1} {\rm Area} \Big( \bigcup_{i=1}^n  C(v_i,t) \Big) 
\]
Thus,  (\ref{conjecture:0721}) is essentially Conjecture~\ref{conj1}. 
For $n=3$, area is circumference and ${\rm Area}(C(v_i,t)) =2 \arccos(t)$. Fix any $W= (w_1, \ldots, w_n)'\in \mathsf{S}$.  If $1/2< t \leq 1$, $C(w_i,t)$ must be disjoint and 
\[
{\rm Area} \Big( \bigcup_{i=1}^3  C(v_i,t) \Big)\leq \sum_{i=1}^3 {\rm Area}(C(v_i,t)) = 6\arccos(t) = {\rm Area} \Big( \bigcup_{i=1}^3  C(w_i,t) \Big)
\]
On the other hand, if $0<t\leq 1/2$, $\bigcup_{i=1}^3  C(w_i,t) $ covers the whole circle, then 
\[
 {\rm Area} \Big( \bigcup_{i=1}^3  C(w_i,t) \Big) = 2\pi \geq {\rm Area} \Big( \bigcup_{i=1}^3  C(v_i,t) \Big)\,.
\]
Thereby, Conjecture~\ref{conj1} and (\ref{conjecture:0721}) hold for $n = 3$. For $n=4$, Conjecture~\ref{conj1} was proved in \cite{landau1966optimality} (see also \cite{balitskiy2017optimality} for some discussions).

 It is important to note that if Conjecture~\ref{conj1} holds, then, by (\ref{conjecture:0721}), $\max_i w_i'u$ first-order stochastically dominates $\max_i v_i'u$ for any $W\in \mathsf S$. Consequently, 
$
\mathbb E [\phi(\max_i w_i'u)] \geq \mathbb E [\phi(\max_i v_i'u)] 
$
 for any increasing measurable function $\phi$ for which the expectations exist. Since $\phi(x) = x$ and $\phi(x) = e^{2rx}$ for any $r>0$ satisfy the conditions,  WS and SMW also follow easily.

 \vspace{0.5em}
\noindent{\bf Proofs of the main theorems}. Theorems~\ref{thm:two-sided}, \ref{thm:hybrid}, and \ref{thm:FWER} are immediate consequences of GCI, since the events involved are all symmetric and convex. Moreover, Conjecture~\ref{conj2} directly implies Theorem~\ref{thm:Slepian}. Therefore, we omit the proofs.

\end{document}